\documentclass[a4paper]{amsart}

\usepackage[utf8]{inputenc}
\usepackage{amsmath, amsfonts, amssymb, amsthm} \parskip1pt
\usepackage{xcolor}
\usepackage{mathrsfs}
\usepackage{amsfonts,graphics,epsfig}
\usepackage{amsmath, amsfonts, amssymb, amsthm, amscd}
\usepackage{hyperref}

\newcommand{\abs}[1]{\ensuremath{\lvert\,#1\,\rvert}}

\newtheorem{theorem}{Theorem}[]
\newtheorem{proposition}[theorem]{Proposition}

\newtheorem{definition}[theorem]{Definition}
\newtheorem{remark}[theorem]{Remark}

\newtheorem{example}[theorem]{Example}

\newtheorem{lemma}[theorem]{Lemma}

\newtheorem{acknowledgement}{Acknowledgement}

\def\R{\mathbb R} \def\Z{\mathbb Z}

\def\N{\mathbb N}
\def\d{\displaystyle}
\def\C{{\mathbb C}} \def\d{{\rm d}}
\def\Q{\mathbb Q}

\def\N{{\mathbb{N}}}

\def\dim{\text {dim}}

\date{}

\def\build#1_#2^#3{\mathrel{\mathop{\kern 0pt#1}\limits_{#2}^{#3}}}

\def\smallsquare{\vbox{\hrule\hbox{\vrule height 1 ex\kern 1
ex\vrule}\hrule}}

\def\N{{\mathbb N}}

\date{}

\begin{document}
\title{On a continued fraction algorithm in finite extensions of $\Q_p$ and its metrical theory}
\author{Manoj Choudhuri and Prashant J. Makadiya}
\address{Department of Basic Sciences, Institute of Infrastructure, Technology, Research and Management, Near Khokhara Circle, Maninagar (East), Ahmedabad 380026, Gujarat, India.}
\email{manojchoudhuri@iitram.ac.in} \email{prashant.makadiya.20pm@iitram.ac.in}
\date{}

\maketitle

\begin{abstract}
We develop a continued fraction algorithm in finite extensions of $\Q_p$ generalising certain algorithms in $\Q_p$, and prove the finiteness property for certain small degree extensions. We also discuss the metrical properties of the associated continued fraction maps for our algorithms using subsequence ergodic theory and moving averages.
\end{abstract}

\vspace{0.2cm}

{\it Mathematics Subject Classification,} Primary: $11$J$70$, $11$J$83$; Secondary: $11$J$61$, $37$A$44$.
   
{\it Keywords:} Continued fractions, finite extensions of $\Q_p$, continued fraction map, exactness, metric theory.

\vspace{0.4cm}

\tableofcontents

\section{Introduction}
Being an indispensable tool in number theory, especially in Diophantine approximation, the study of continued fractions has attracted many mathematicians over the years. The simple or classical continued fraction expansion of a real number $\alpha$ is an expression of the form 
\begin{align}\label{SCFE}
\alpha=a_0+\frac{1}{\displaystyle{a_1+\frac{1}{\displaystyle{a_2+\frac{1}{\displaystyle{a_3+_{\ddots} }}}}}},   
\end{align}
which is also written as $\alpha=[a_0;a_1,a_2,\dots]$ with $a_i$'s being natural numbers and $a_i>0$ for $i\geq 1$ (see \cite{HW} or \cite{Kh} for more details). Here, $a_i$'s are called the partial quotients of the continued fraction expansion of $\alpha$. If $\frac{A_n}{B_n}=[a_0;a_1,\dots,a_n]$, then the rational numbers $\frac{A_n}{B_n}$ converges to $\alpha$, and $\frac{A_n}{B_n}$ is called the $n$th convergent to the continued fraction of $\alpha$. The classical continued fraction for real numbers has nice arithmetical properties such as rational numbers have finite continued fraction expansion; convergents are the best approximants among other rational numbers; quadratic irrationals have periodic continued fraction expansion expansions and vice versa, this fact is known as Lagrange's theorem (see \cite{Kh} for more details). The reader may look at \cite{KU1} for various real continued fractions apart from the classical (simple) one.
The starting of continued fraction theory for complex numbers goes back to $1887$ when A. Hurwitz (\cite{Hur87}) described the nearest integer continued fraction algorithm in the field of complex numbers, the partial quotients being elements of the ring of Gaussian integers. He also proved a version of lagrange's theorem as well. See \cite{D2014}, \cite{D2015} for more recent developments and a general approach to continued fraction theory in this setup. 

It is quite natural to study continued fractions in the non-Archimedean setup as well. The reader is referred to \cite{WSch1} for a comprehensive introduction to the theory of continued fraction and its relation to Diophantine approximation in positive characteristics. See also the survey article \cite{L1} by Lasjaunias. For continued fraction in the field of Laurent series in one indeterminate over a finite field $\mathbb{F}_q$, viz. $\mathbb{F}_q((X^{-1}))$, there is a natural choice of a set for the set of partial quotients, viz. the polynomial ring $\mathbb{F}_q[X]$. The continued fraction in this setup is very well-behaved. For example, any element in $\mathbb{F}_q(X)$ has a finite continued fraction expansion, the convergents (which naturally belong to $\mathbb{F}_q(X)$) provide the best approximation, in fact it is true that if $\alpha\in\mathbb{F}_q((X^{-1}))$ and $\frac{P}{Q}\in \mathbb{F}_q(X)$ is such that $\left|\alpha-\frac{P}{Q}\right|<\frac{1}{|Q|^2}$, then $\frac{P}{Q}$ is a convergent from the continued fraction expansion of $\alpha$. A version of Lagrange's theorem is true as well in this setup, see \cite{WSch1} or \cite{L1} for more details.

In $1940$, Mahler (\cite{Ma1}) initiated the study of continued fractions in the field of p-adic numbers. There are mainly two types of continued fractions in the field of $p$-adic numbers, one was introduced by Schneider (\cite{Sc1}) in $1968$, and the other was introduced by Ruban \cite{RA1} in $1970$. Rational numbers need not have finite continued fraction expansion with respect to these algorithms. See \cite{Bu} for rational numbers having infinite expansion with respect to Schneider's algorithm. In fact, Wang \cite{WL1} and Laohakosol \cite{Lao85} independently showed that a $p$-adic number $\alpha$ is rational if and only if the Ruban continued fraction expansion of $\alpha$ is either finite or periodic. In $1978$, Browkin (\cite{BJ1}) modified Ruban's algorithm and proved that every rational number has a finite continued fraction expansion. Another desirable arithmetical property of any continued fraction is periodic expansion (the periodicity property) of quadratic irrational which is known as Lagrange's theorem in the case of real numbers. In \cite{BJ2}, Browkin modified his algorithm further and showed that $\sqrt{m}$ has periodic continued fraction expansion for certain positive integers for $p=5$. Though, the same is not true for larger values of $p$. So, Lagrange's theorem is not true in this setup. See also references cited there in \cite{BJ2} for various work related to periodicity prior to Browkin's (\cite{BJ2}) work in $2000$. Many research works have been done in recent times in which people have presented many modified algorithms to achieve the periodicity and other desirable properties of continued fractions in the field of $p$-adic numbers. See for example \cite{CVZ19}, \cite{NLL1},\cite{CMT2023},\cite{MRS2023}, \cite{MR2024}, and the references cited there in. See also the survey article by Romeo (\cite{Ro}) for a comprehensive history of the development of the continued fraction theory in the field of $p$-adic numbers.

It is quite natural to consider continued fraction in finite extensions of $\Q_p$ which is left-out in the above discussion while considering continued fractions in all locally compact fields.
In this article, we consider canonical extensions of  the algorithms of Ruban and Browkin for $\Q_p$ in its finite extensions. Given any finite (necessarily simple) extension $K$ of $\Q_p$, we consider this extension in two steps, viz.  $K=L(\beta)$ with $K/L$ a totally ramified extension, and $L=\Q_p(\gamma)$ with $L/\Q_p$ an unramified extension (see next section for more details). In our algorithms, the partial quotients are elements from the set $\Z\left[\frac{1}{p}\right][\gamma,\beta]$. 
We show that any $\alpha \in K$ has a unique continued fraction expansion, and given any sequence of partial quotients $\{c_i\}_{i \geq 0}$, $[c_0,c_1,\ldots,c_n]$ converges to an element $\alpha$ of $K$. For a few small degree unramified extensions of the form $\mathbb{Q}_p(\gamma)$, we show that any element of $\mathbb{Q}(\gamma)$ has a finite continued fraction expansion.

In this article, we are also going to discuss the metrical theory of the associated continued fraction map. In the case of classical continued fraction for real numbers the Gauss map or the continued fraction map is defined as
\[
T: (0,1) \rightarrow (0,1)
\]
\[
T(x)=\left\{\frac{1}{x}\right\},
\]
where $\left\{\frac{1}{x}\right\}$ denotes the fractional part of $\frac{1}{x}$.
It is well known that $T$ is ergodic with respect to the Gauss measure (see \cite{EW}). For the ergodicity of the continued fraction map for a more general class of continued fraction, see \cite{KU2}.
For complex continued fraction, the ergodicity of the maps associated with the nearest integer complex continued fractions over imaginary quadratic fields is discussed in a recent paper of Nakada et al. \cite{ENN}. They, in fact, showed that the continued fraction map is exact (see Section $4$ for definition). See also some references in \cite{ENN} for some earlier works related to the metrical theory of complex continued fractions.

In the non-Archimedean settings, Berthe and Nakada \cite{BN1} proved the ergodicity of the continued fraction map in positive characteristic, and as an application they obtained various metrical results regarding the averages of partial quotients, average growth rates of the denominators of the convergents, etc. In \cite{LN14}, Lertchoosakul and Nair proved the exactness of the continued fraction map using which they could consider more general averages concerning the partial quotients and the growth rate of the denominator of the convergents. The quantitative version of the metric theory of the continued fraction map in this setup was considered by the same authors in a subsequent paper \cite{LN18}. The reader is referred to \cite{V1964} for quantitative metrical results concerning real continued fraction.

In this article, we discuss metrical theory of continued fractions in finite extension of $\Q_p$. We show that the associated continued fraction map is Haar measure preserving and exact. Then we obtain various metrical results analogous to the results of \cite{LN14} concerning asymptotic behaviour of various quantities related to partial quotients, denominator of the convergents, etc. In these results, general averages using subsequence ergodic theory and moving averages are considered as done in \cite{LN14}.  

\section{Preliminary}
For a prime number $p$, the field of $p$-adic numbers $\Q_{p}$ is the set of all Laurent series in $p$ of the form
  \[
    \alpha=\sum\limits_{j\geq n_0 }a_jp^{j}, \text{ where } a_j \in \{ 0, 1,\ldots,p-1 \} \text{ and } n_0 \in \mathbb{Z}.
  \]
  The $p$-adic valuation $\mathfrak{v}_p$ on $\Q_{p}$ is defined as follows: if $\alpha=\sum_{j\geq n_0 }a_jp^{j}$, then
  \[
    \mathfrak{v}_p(\alpha) := \inf\,\{\,j \in \mathbb{Z}\ :\ a_j \neq 0\,\}.
\]
  Then the $p$-adic absolute value of $\alpha$ is given by 
  \[
  |\alpha|_p := p^{-\mathfrak{v}_p(\alpha)}
  \]
  when $\alpha\neq 0$, and $|0|_p=0$. The field of $p$-adic numbers is the completion of $\Q$ with respect to this absolute value.
 Let $K = \Q_p(\xi)$ be a finite extension of $\Q_p$ of degree $m$, i.e., $[\,K\,:\,\Q_p\,] = m$. We may then take $\mathfrak{B}=\{ 1, \xi ,\ldots, \xi^{m-1} \}$ as a convenient vector space basis for $K$ over the f{i}eld $\Q_p$. Otherwise said, any element $\mathfrak{b} \in K$ can be written uniquely as
  \[
    \mathfrak{b} =b_{0} + b_{1}\xi + \cdots + b_{m-1}\xi^{m-1}, \text{ where } b_j \in \mathbb{Q}_{p} \text{ for all } j.
  \]
Since every finite extension is an algebraic extension, we have for every $\mathfrak{b} \in K$ that there is some monic irreducible polynomial
  \[
    g(x) = x^{n} + B_1 x^{n-1} + \cdots + B_{n-1} x + B_n
  \]
  of degree at most $m$ and coefficients $B_j \in \mathbb{Q}_{p}$ such that $g(\mathfrak{b})=0$ in $K$. The \emph{norm map} for the f{i}nite f{i}eld extension $K / \Q_p$ is then def{i}ned as
  \[
    N_{K / \mathbb{Q}_{p}} (\mathfrak{b}) := (-1)^{n} B_{n}.
  \]
  Our absolute value $| \cdot |_p$ on $\Q_p$ extends uniquely to $K$ in the following manner (see~\cite{NK} for details):
  \[
    |\mathfrak{b}| := \big| N_{K / \mathbb{Q}_{p}} (\mathfrak{b}) \big|_p^{\frac{1}{n}},\quad \mathfrak{b} \in K.
  \]
  Let us choose an element $\pi \in K$ of maximal absolute value smaller than 1, say $0 < \theta := |\pi| < 1$. 
 % Note that the existence of such a $\pi$ is guaranteed as the value group $\mathfrak{v}_p (K^*)$ is a discrete subset of $\mathbb{R}$.
  Def{i}ne
  \[
    \mathcal{O}_K := \{\,x\in K\ :\ |x| \leq 1\,\},\quad\mathfrak{m}_K := \{\,x\in K\ :\ |x| < 1\,\}
  \]
  and
  \[
    \mathcal{O}_K^* := \{\,x \in K\ :\ |x| = 1\,\}.
  \]
  We have $\pi \mathcal{O}_K = \mathfrak{m}_K$, and the residue field $\overline{k} = \mathcal{O}_K/\mathfrak{m}_K$ is a finite extension of $\mathbb{F}_p$.
  \begin{definition}
    The \textbf{residue degree} of the finite extension $K$ of $\mathbb{Q}_p$ is the positive integer $f=[\overline{k}:\mathbb{F}_p]=\dim_{\mathbb{F}_p}(\overline{k})$, where $\overline{k}$ is the residue field of $K$. A finite extension $K$ of $\mathbb{Q}_p$ is said to be \textbf{totally ramified} if $f=1$. 
  \end{definition}
  We also have that $\overline{k}=\mathbb{F}_q$, where $q = \#\,( \overline{k} ) = p^f$. This is because upto isomorphism there is exactly one f{i}nite f{i}eld having $q$ elements.
  \begin{definition}
    The \textbf{ramification index} of $K / \mathbb{Q}_p$ equals $e = [\,|K^*|\,:\,|\mathbb{Q}^{*}_{p}|\,] = \#( |K^*| / p^{\mathbb{Z}} )$. A finite extension $K$ of $\mathbb{Q}_p$ is called \textbf{unramified} if $e=1$. 
  \end{definition}
  
  For the \emph{uniformizer} $\pi \in \mathfrak{m}_K$, we have $|\pi|^e = |p|$ thereby giving us
  \[
    |\pi| = p^{- 1 / e}.
  \]
  Any element $\alpha \in K$ can be represented as $\alpha = u\pi^n$ for some suitable $u\in \mathcal{O}^{*}_K$ and $n \in \mathbb{Z}$. Then,
  \[
    |\alpha| = |\pi|^n = p^{- n / e}.
  \]
  The integers $e$ and $f$ given above sat{i}sfy $ef = m$, where $m$ is the degree of the extension. We recall that ~\cite[Corollary 4-26]{RV1} there exists an unrami{f}ied subextension $L / \Q_p$ of degree $f$ such that $K / L$ is a totally rami{f}ied extension of degree $e$. Also, there exists a $\gamma \in L$ such that $L=\mathbb{Q}_p(\gamma)$ with $|\gamma|=1$~\cite[\S\,III.3]{NK}. (To boot, we may and do take $\gamma$ to be some primi{t}ive $(p^f-1)$-th root of unity)
  \begin{lemma}\label{L:beta}
    Let $K / L$ be as above. Then, there exists some $\beta \in K$ such that $K=L(\beta)$ and $|\beta| = p^{1 / e}$.
  \end{lemma}
  \begin{proof}
    We know that the value groups of $L^*$ and $K^*$ are $\Z$ and $( 1 / e ) \Z$, respect{i}vely. Let us, therefore, choose some $\beta \in K \setminus L$ such that $| \beta | = p^{1/e}$. As $K / L$ is a f{i}nite extension, every element of $K$ is algebraic over the f{i}eld $L$. In part{i}cular, our chosen element $\beta$ sat{i}sf{i}es some minimal monic polynomial
    \[
      h (x) = x^n + b_{n - 1}x^{n - 1} + \cdots + b_0,\ b_j \in L
    \]
    and $n \leq e$.
    % As discussed before, we may assume that $K = L(\xi)$ for some $\xi \in K$. This element $\xi$ can in fact be chosen to have absolute value $\geq 1$. If $|\xi| = 1$, let us take $\beta = \xi$ and we are done. If not, suppose that $|\xi| = p^{a / e}$ where $a$ is some posi{t}ive integer (as $|\xi|\geq 1$).\\[-0.2cm]
    % When $a$ is divisible by $e$, we may multiply $\xi$ by $p^{a / e } \in \Q_p \subset L$ and obtain a new generator $p^{a / e}\xi$ for the f{i}eld extension $K / L$. Hence, the only case that remains to be discussed is when $a$ is not divisible by $e$. Without loss of generality, we may assume $0 < a < e$ by repeat{i}ng the last argument.\\[-0.2cm]
    % 
    % Suppose $\gcd ( a, e ) = d$. Then, there exist integers $b$ and $t$ such that $\frac{ba}{e}=\frac{d}{e} + t$ and $\delta$ in $K$ such that $|\delta|=p^{-t+\frac{1-d}{e}}$. Let $\beta=\delta \xi^q$. Then $|\beta|=p^{\frac{1}{e}}$.\\
    Now, we will like to show that $K = L(\beta)$. It is equivalent to establishing that the minimal polynomial of $\beta$ has degree $e$. Suppose $n < e$.\\[-0.2cm]
    
% \par Suppose the minimal polynomial of $\beta$ is given by $$b_0+b_1\beta +\cdots +b_{n-1}\beta^{n-1}+\beta^n; \ b_i \in L, \ n<e.$$
    We have $|\beta|=|b_0|^{1 / n}$ by the unique extension of the non-archimedean absolute value to $K$. Here, $|b_0| = p^s$ for some $s\in \mathbb{Z}$ implying that $|\beta|=p^{s / n} = p^{1 / e}$. This is possible if{f} $n = se$ but $0 < n < e$, a contradiction. Thus, there exists a $\beta \in K$ such that $K=L(\beta)$ with $|\beta|=p^{1 / e}$.
  \end{proof}
  \noindent Every $\alpha \in K=L(\beta)=\mathbb{Q}_p(\gamma)(\beta) = \Q_p ( \beta, \gamma )$ can then be wri{t}ten as
\begin{align}\label{eq11}
\alpha = \sum_{i = 0}^{e - 1} \sum_{j = 0}^{f - 1} b_{i, j} \gamma^j \beta^i,\ b_{i, j} \in \mathbb{Q}_p.
\end{align}
Now let $X_1$ and $X_2$ be the sets inside $\Q_p$ defined by: 
\begin{align}\label{eq12}
X_1 = \left\{\, \sum \limits_{j=0}^{k}\frac{a_j}{p^j}:\ k\in \N\cup \{0\}\ \text{and}\ a_j \in \{ 0, 1, \ldots, p - 1 \} \ \text{for} \ 0 \leq j \leq k \right\},
\end{align}
and 
\begin{align}\label{eq12*} 
X_2 = \left\{\,\sum \limits_{j=0}^{k}\frac{a_j}{p^j}:\ k\in \N\cup \{0\}\ \text{and}\ a_j \in \Big\{ -\frac{p-1}{2}, \ldots, \frac{p-1}{2} \Big\} \ \text{for} \ 0 \leq j \leq k \right\}.
\end{align}
Note that the partial quotients for Ruban's $p$-adic continued fraction are elements of $X_1$, whereas, the partial quotients for Browkin's algorithm are elements of $X_2$. In a moment, we define a set $Z$, the elements of which will be used as partial quotients for the continued fraction algorithm developed in this article. We may either use $X_1$ or $X_2$ while defining the set $Z$. In the first case, we get a generalization of Ruban's algorithm in finite extensions of $\Q_p$, whereas, we get a generalization of Browkin's algorithm in the second case. From now on, we use the notation $X$ for both the sets $X_1$ and $X_2$ with the understanding that whenever we use $X$, the discussion applies to both $X_1$ and $X_2$. 
Any two dist{i}nct numbers in $X$ will have $p$-adic distance at least one giving us that it is a \textit{$1$-uniformly discrete set.} Furthermore, every non-zero element has an absolute value at least one. Generalizing this observat{i}on for all f{i}nite $p$-adic extensions, we have the following.
  \begin{lemma}\label{L:1ud}
    The set
    \[
      Z := \left\{\,\sum_{i=0}^{e-1} \sum_{j=0}^{f-1} b_{i,j}\gamma^j \beta^i\ :\ \ b_{i,j} \in X\,\right\}
    \]
    is $1$-uniformly discrete. In part{i}cular, every non-zero element in $Z$ has an absolute value at least one.
  \end{lemma}
  \begin{proof}
    Let $z_1, z_2 \in Z$ with $z_1 \neq z_2$. This happens if{f}
    \[
      z_1 - z_2 = \sum_{i=0}^{e-1} \sum_{j=0}^{f-1}\,\widetilde{b}_{i,j}\gamma^j \beta^i,
    \]
    where at least one of the coef{f}icients $\widetilde{b}_{i, j}$'s (say $\widetilde{b}_{k, \ell}$) is a non-zero element from the set $X \cup ( -X )$. Consider
    \begin{equation}\label{E:yinL}
      y = \sum_{j = 0}^{f - 1} \widetilde{b}_{k, j} \gamma^j \in L.
    \end{equation}
    Assume $|y| < 1$. Since $\gamma$ is a primi{t}ive $f$-th root of unity, it is plain that $\gamma$ belongs to $\mathcal{O}_L$. Without loss of generality, we may assume that $\widetilde{b}_{k,f-1}$ has the maximum absolute value in the representat{i}on of $y$ given by~\eqref{E:yinL}. This is because by choosing $m$ such that
    \[
      \abs{\widetilde{b}_{k, m}} = \max_{0\,\leq\,j\,\leq\,f - 1} \big\{\,\abs{\widetilde{b}_{k, j}}\,\big\}
    \]
    and replacing $y$ with $\gamma^{f - 1 - m}y$, we can ensure that the coef{f}icient of $\gamma^{f-1}$ has maximum absolute value amongst all the $\widetilde{b}_{k,j}$'s. It then follows that
    \begin{equation}\label{E:gammay}
      \gamma^{f-1} = ( \widetilde{b}_{k, f - 1} )^{-1} \left( y-\sum_{j = 0}^{f - 2} \widetilde{b}_{k, j} \gamma^j \right)
    \end{equation}
    On reducing the above equat{i}on modulo $p\mathcal{O}_L$, we get that 
    \begin{equation}\label{eqn3}
      \gamma^{f-1}=\sum_{j = 0}^{f - 2} a_j \gamma^j
    \end{equation}
    where $a_j$'s are elements of the residue field $\mathbb{Z}_p / p\mathbb{Z}_p$. This leads to a contradiction as the degree of the residue f{i}eld $\mathcal{O}_L/p\mathcal{O}_L$ over $Z_p/pZ_p$ is $f$. Therefore, $\abs{y} \geq 1$. % 
  We will in fact have that
  \[
    z_1 - z_2 = \sum_{i=0}^{e-1}y_i \beta^{i}
  \]
  with $|y_i| = p^{n_i}$ for some $n_i \in \N \cup \{ 0 \}$ or $\abs{y_i} = 0$ for some $i$ while $|\beta|=p^{1 / e}$. This implies that $\abs{y_{i_1} \beta^{i_1}} \neq \abs{y_{i_2} \beta^{i_2}}$ for any pair of indices $0 \leq i_1 < i_2 < e$ with at least one $y_{i_{1}}$ or $y_{i_{2}}$ non-zero and there exists $k\in \{ 0,\ldots,e-1 \}$ such that
  \[
    \max _{0 \leq i \leq e-1} \{ |y_i\beta^{i}| \} =|y_{k}\beta^k| \geq 1.
  \]
  Therefore, 
  \begin{align*}
      |z_1-z_2|&=\left| \sum_{i=0}^{e-1}y_i \beta^{i} \right| \\&=\max _{0 \leq i \leq e-1} \{ |y_i\beta^{i}| \} \\&=|y_{k}\beta^k|\\& \geq 1.
  \end{align*}
  Hence, $Z$ is a $1$-uniformly discrete set.
\end{proof}
The metric balls in $K$ will have radius $\frac{s}{e}$ for $s\in\Z$. More precisely, for $\alpha \in K$ and $s \in \Z$, let
\[
    B\,(\,\alpha,\,p^{\frac{s}{e}}\,)\ =\ \{\,x \in K\ :\ | x - \alpha | < p^{\frac{s}{e}}\,\}
 \]
be the ball around $\alpha$ of radius $p^{\frac{s}{e}}$. Let $\mu$ denote the Haar measure on the local f{i}eld $K$ (see Chapter $4$ of \cite{RV1} for existence of Haar measure) and it is normalized in such a way that $\mu\,\big(\,B ( 0, 1 )\,\big) = 1$. Note that
  \[
    B ( 0, p^s )\ =\ p^{-s} B ( 0, 1 )\ =\ \{\,p^{-s}x\ :\ x \in B ( 0, 1 )\,\}
  \]
and, therefore,
  \begin{align}
    \mu\,\big(\,B ( 0, p^s )\,\big) &= \mu \big(\,p^{-s}B(0,1)\,\big) =\!\mod_K\left(\frac{1}{p^s}\right) \mu (B(0,1))\label{E:scale}\\
    &=\!\mod_{\Q_p} \left( \frac{1}{p^s} \right)^m\ =\ p^{sm}\notag
  \end{align}
  by~\cite[\text{Proposition 4-13}]{RV1}. 
 Next two technical lemmas are useful for computing measure of various metric balls inside $K$.
  \begin{lemma}\label{lem5}
    For $s\in \Z$, the ball $B ( 0, p^s ) \subset K$ is the same as the set
    \[
      A := \big\{\,x \in K\ :\ x = \sum\limits_{i=0}^{e-1} \sum \limits_{j=0}^{f-1} b_{i,j} \gamma^j \beta^i,\ b_{i,j} \in p^{1-s} \Z_p\,\big\}.
    \]
  \end{lemma}
  \begin{proof}
    F{i}rst suppose $x \in A$. Then, $x\ =\ \sum_{i=0}^{e-1} \sum_{j=0}^{f-1} b_{i,j} \gamma^j \beta^i$ where $b_{i,j} \in p^{1-s} \Z_p$ for all $i, j$. Since $|b_{i,j}| \leq p^{s-1},\ |\gamma| = 1$ and $|\beta|=p^{\frac{1}{e}}$, we obtain $|x| \leq p^{s-1+\frac{e-1}{e}}<p^s$. Otherwise said, $A \subseteq \{ x \in K : |x|<p^s \} = B ( 0, p^s )$.\\[-0.2cm]

    Conversely, let $x \notin A$ so that $x = \sum_{i=0}^{e-1} \sum_{j=0}^{f-1} b_{i,j} \gamma^j \beta^i$, where $|b_{i,j}| \geq p^s$ for at least one pair $( i, j )$. For each $0 \leq i < e$ and $ 0 \leq j < f$, we may write
    \[
      b_{i , j} = \{ b_{i, j} \} + [ b_{i, j} ] \in \Q_p
    \]
    with all terms of the form $\sum_{l=0}^{k}a_{-s-l}p^{-s-l}$ contained in $\{ b_{i,j} \}$, $[ b_{i, j}] \in p^{1-s} \Z_p$ and $k \in \N  \cup \{ 0 \}$ is such that $a_{-s-l}=0$ for all $l>k$.
    Thus, $x = x_1 + x_2$ where $p^s x_1$ belongs to the set $Z$ introduced in Lemma~\ref{L:1ud} while $x_2 \in A$. This implies that $\abs{x_1} \geq p^s$ and $\abs{x_2} < p^s$. All in all, $\abs{x} \geq p^s$ or equivalently, $x \notin B ( 0, p^s )$.
  \end{proof}
\begin{lemma}\label{lem6}
    Let $s \in \Z$ and $a \in \{ 1,\ldots,e-1 \}$. Then the ball $B(0,p^{s-1+\frac{e-a}{e}}) \subset K$ is the same as the set 
    \[
    A := \big\{ \,x \in K: \,x = \sum\limits_{i=0}^{e-1} \sum\limits_{j=0}^{f-1} b_{i,j} \gamma^j \beta^i, \begin{matrix} 
     \ b_{i,j} \in p^{1-s} \Z_p \ \text{for} \ i < e - a \ \text{and} \\
    b_{i , j} \in p^{2-s} \Z_p \ \text{for} \ e-a \leq i <e
    \end{matrix}   \big \}
    \]
\end{lemma}
 \begin{proof}
     First we show that $A \subseteq B(0, p^{s-1+\frac{e-a}{e}})$. Let $x \in A$. Then $x = \sum\limits_{i=0}^{e-1} \sum\limits_{j=0}^{f-1} b_{i,j} \gamma^j \beta^i$, where 
     $b_{i,j} \in p^{1-s} \Z_p$ for $i < e - a$, and $b_{i , j} \in p^{2-s} \Z_p$ for $e-a \leq i <e$. Note that $|b_{i,j}| \leq p^{s-1}$ for all $0 \leq j <f$, $0 \leq i <e-a$ and $|b_{i,j}| \leq p^{s-2}$ for all $0 \leq j <f$, $e-a \leq i <e$. Also, $|\gamma|=1$ and $|\beta|=p^{\frac{1}{e}}$, we get $|x| \leq p^{s-1+\frac{e-a-1}{e}}<p^{s-1+\frac{e-a}{e}}$. This shows that $x \in B(0,p^{s-1+\frac{e-a}{e}})$.

To prove the converse, suppose $x \notin A$. Then there exists a pair $(i,j)$ such that 
    \[
      b_{i,j} \notin p^{1-s} \Z_{p} \ \text{for} \ 0 \leq i <e-a \ \text{or} \  b_{i,j} \notin p^{2-s} \Z_{p} \ \text{for} \ e-a \leq i <e.
    \]
    If $b_{i,j} \notin p^{1-s} \Z_{p}$, then $|b_{i,j}| \geq p^{s}$. So, $|x| \geq p^{s} \nless p^{s-1+\frac{e-a}{e}}$. And if $b_{i,j} \notin p^{2-s} \Z_{p}$, then $|b_{i,j}| \geq p^{s-1}$. So, $|x| \geq p^{s-1+\frac{e-a}{e}} \nless p^{s-1+\frac{e-a}{e}}$. Thus, $x \notin B(0,p^{s-1+\frac{e-a}{e}})$.
 \end{proof} 
 
Now we are in a position to calculate the measure of any ball around zero inside $K$.
\begin{proposition}\label{measure}
The measure of the ball $B \big( 0,\ p^{s+ \frac{i}{e}} \big)$ equals $p^{sm + f{i}}$.
\end{proposition}
%we calculate $\mu(B(0,p^{s+\frac{i}{n}}))$, where $s \in \Z$ and $i \in \{1,\ldots,e-1 \}$. 
\begin{proof}
Fix some $a \in \{ 1,\ldots,e-1 \}$. By (\ref{E:scale}) we know that
  \begin{align*}
    p^{sm}\ &=\ \mu\,\big(\,B(0,p^s)\,\big)\ =\ \mu\,\big(\,\{\,x \in K\ :\ \abs{x} < p^s\,\}\,\big)\\
    &=\ \mu\,\big(\,\{\,x \in K\ :\ x =\sum \limits_{i=0}^{e-1} \sum \limits_{j=0}^{f-1} b_{i,j} \gamma^j \beta^i,\ b_{i,j} \in p^{1-s}\Z_p\,\}\,\big) \ \ (\text{by Lemma} \ \ref{lem5} ).
  \end{align*}
  Since $p^{1-s} \Z_p / p^{2-s} \Z_p \simeq ( \{ 0,p^{1-s},\ldots,(p-1)p^{1-s} \} , +,*)$, where
  \[
  t_jp^{1-s}+t_{j'}p^{1-s}:= (t_j +_{p}t_{j'})p^{1-s} \ \text{and} \ t_jp^{1-s}*t_{j'}p^{1-s}:= (t_j *_{p}t_{j'})p^{1-s}
  \]
  for all $t_j, t_{j'} \in \{ 0,1,\ldots,p-1 \}$, we have a disjoint union decomposition
  \[
  p^{1-s} \Z_p = p^{2-s} \Z_p \sqcup (p^{1-s}+p^{2-s} \Z_p) \sqcup \ldots \sqcup ((p-1)p^{1-s} +p^{2-s} \Z_p).
  \]
  Therefore,
  \begin{align*}
      p^{sm} \ &=\ \mu \left( \bigsqcup \big\{ x =\sum\limits_{i=0}^{e-1} \sum \limits_{j=0}^{f-1} b_{i,j} \gamma^j \beta^i, \begin{matrix} b_{i,j} \in p^{1-s} \Z_p \ \text{for} \ i < e - a \ \text{and} \\
     b_{i , j} \in t_j p^{1-s} +p^{2-s} \Z_p \ \text{for} \ e-a \leq i <e
    \end{matrix}   \big\} \right),
  \end{align*}
  where the disjoint union is taken over all possible combinat{i}ons\\ $t_j \in \{\,0, 1, \ldots, p - 1\,\}$. By the translation invariant property of $\mu$,
  \[
  \mu (t_j p^{1-s} +p^{2-s} \Z_p)=\mu (p^{2-s} \Z_p).
  \]
  So,
  \begin{align*}
    p^{sm}\ &=\ p^{fa} \mu \left( \big\{\,x = \sum\limits_{i=0}^{e-1} \sum\limits_{j=0}^{f-1} b_{i,j} \gamma^j \beta^i, \begin{matrix} 
     \ b_{i,j} \in p^{1-s} \Z_p \ \text{for} \ i < e - a \ \text{and} \\
    b_{i , j} \in p^{2-s} \Z_p \ \text{for} \ e-a \leq i <e
    \end{matrix}   \big\} \right)\\
    &=\ p^{fa} \mu\,\big(\,B (\,0,\,p^{s-1+\frac{e-a}{e} }\,)\,\big) \ \ (\text{by Lemma} \ \ref{lem6}).
\end{align*}
It follows that $\mu\,\big(\,B\,( 0, p^{s-1+\frac{e-a}{e}} )\,\big) = p^{sm - fa}$. In general, 
\[
\mu \big( B ( 0, p^{s+ \frac{i}{e}} ) \big)=p^{sm + f{i}}
\]
for all $s \in \Z$ and $i \in \{\,0, \ldots, e - 1\,\}$.
\end{proof}
As $Z$ is a uniformly discrete set, we can count the number of elements inside a set of elements with a fixed absolute value. This counting will be useful in some of the subsequent sections. Let $Y$ be the set given by  
\[
Y := \left \{ \sum_{j=0}^{f-1}b_{j} \gamma^{j} :b_{j} \in X \right \} .
\]
For $y \in Y$, $|y| \leq p^t$ for some $t \in \N \cup \{ 0 \}$ if and only if $|b_j|<p^{t+1}$ for all $j=0,\ldots,f-1$. Then it follows that for each $n \geq 1$, 
\begin{align*}
    \# \{ y \in Y:|y|=p^n \} &= \# \{ y \in Y:|y| \leq p^n \} - \# \{ y \in Y:|y|<p^n \} \\
    &=\big( p^{n+1} \big)^f - \big( p^{n} \big)^f \\
    &=p^{fn}(p^f-1).
\end{align*}
Let us denote by $Z^*$ the set of all $c \in Z$ such that $|c|>1$. Also let $s$ be a positive integer and $a \in \{ 0,\ldots,e-1 \}$. Writing an element $c\in Z^*$ as
\[
c=\sum_{i=0}^{e-1} \sum_{j=0}^{f-1} b_{i,j}\gamma^j \beta^i\ :\ \ b_{i,j} \in X,
\]
note that if $|b_{i,j}| \geq p^{s+1}$ for some $0 \leq i \leq a \ \& \ 0 \leq j<f$ or $|b_{i,j}| \geq p^{s}$ for some $a<i<e \ \& \ 0 \leq j<f$, then the set $\{ c \in Z^{*}:|c| \leq p^{s+\frac{a}{e}} \} =\phi $. Then it follows that 
\begin{align*}
\# \left \{ c \in Z^{*}:|c|=p^{s+\frac{a}{e}} \right \} &=\# \left \{ c \in Z^{*}:|c| \leq p^{s+\frac{a}{e}} \right \} - \# \left \{ c \in Z^{*}:|c| \leq p^{s+\frac{a-1}{e}} \right \} \\&= \big(p^{s+1}\big)^{f(a+1)} \big(p^{s}\big)^{f(e-(a+1))} \, - \, \big(p^{s+1}\big)^{af}\big(p^{s}\big)^{f(e-a)} \\&=p^{f(a+es)}\big(p^f-1 \big).    
\end{align*}
Now, let $a \in \{ 1,\ldots,e-1 \}$. By a similar argument, we also obtain
\[
\# \{ c \in Z^{*}:|c|=p^{\frac{a}{e}} \} =p^{af}(p^f-1).
\]
Hence, for all $n \in \N$,
\begin{equation}\label{Z*:coun}
    \# \big \{ c \in Z^{*}:|c|=p^{\frac{n}{e}} \big \}=p^{fn}\big( p^f-1 \big).
\end{equation}

\section{Continued fraction algorithm and finiteness property}
Now we describe a continued fraction algorithm for elements in any finite extension $K$ of $\Q_p$. This algorithm generalizes Ruban's algorithm of $\Q_p$ (for $X=X_1$), as well as Browkin's algorithm (for $X=X_2$). Our algorithm is a very natural extension of Ruban's and Browkin's algorithm with partial quotients coming from the set $Z$.

The $p$-adic floor function for Ruban's algorithm is a function from $\Q_p$ to $X$ defined as follows: for $\alpha=\sum \limits_{j \geq n_0}a_jp^j \in \Q_p$ with $a_j \in \{ 0,1,\ldots,p-1 \}$ for $X=X_1$ and $a_j \in \{ -\frac{p-1}{2},\ldots,0,\ldots,\frac{p-1}{2} \}$ for $X=X_2$,
\begin{equation*}
\ \lfloor \alpha \rfloor _p=\left\{ \begin{array}{rcl} \displaystyle{\sum \limits_{j=n_{0}}^{0}a_{j}p^{j}}, & \mbox{if} & \mathfrak{v}_{p}(\alpha) \leq 0 \ (\text{or} \ |\alpha|_{p}\geq 1) \\ 0 \ \ \ , & \mbox{} & \text{otherwise.}
\end{array}\right.  
\end{equation*}
Using this floor function, we define a floor function on $K$ which is a function from $K$ to $Z$, as follows:
\[
\text{For} \ \alpha = \sum_{i = 0}^{e - 1} \sum_{j = 0}^{f - 1} b_{i, j} \gamma^j \beta^i,\ b_{i, j} \in \mathbb{Q}_p,
\]
\begin{equation}\label{E:floor}
\lfloor \alpha \rfloor=\sum_{i = 0}^{e - 1} \sum_{j = 0}^{f - 1} \lfloor b_{i, j} \rfloor_p\  \gamma^j \beta^i .    
\end{equation}
It is easy to see that 
\begin{equation}\label{**}
|\alpha - \lfloor \alpha \rfloor| \leq \frac{1}{p^{\frac{1}{e}}}<1
\end{equation}
for any $\alpha \in K$.

Following the existing literature, we call an expression of the form 
\[
\displaystyle{c_0+\frac{1}{\displaystyle{c_1+\frac{1}{\ddots}}}}
\]
with $c_j \in Z$ and $c_j \in Z^{*}$ for $j \geq 1$, a continued fraction which is also written as $[c_0;c_1,\ldots]$. It is a finite continued fraction if the sequence $(c_j)$ is a finite one, otherwise, it is an infinite continued fraction. We call $c_j$'s the partial quotients of the continued fraction. We write,
\[
[c_0;c_1,\ldots,c_n]=\frac{s_n}{t_n},
\]
with $s_n$, $t_n \in \Q(\gamma , \beta)$, and call it the $n$th convergent of the continued fraction $[c_0;c_1,\ldots]$. It is easy to see that the sequence $(s_n)$ and $(t_n)$ satisfy the following recurrence relations:
\[
s_n=c_ns_{n-1}+s_{n-2} \ \text{and} \ t_n=c_nt_{n-1}+t_{n-2}, \ n \geq 2,
\]
with $s_0=c_0$, $t_0=1$, $s_1=c_0c_1+1$, $t_1=c_1$. The numerator and denominator of the convergents also satisfy
\begin{equation}\label{*}
t_ns_{n-1}-s_nt_{n-1}=(-1)^{n}, \ n \geq 0.
\end{equation}
Now we discuss the convergence properties of the continued fraction in our setup.
\begin{lemma}\label{+}
Let $c_{0},c_{1},\ldots \in Z$ with $c_j \in Z^{*}$ for $j \geq 1$. Then the sequence of convergents $\frac{s_n}{t_n}=[c_0;c_1,\dots,c_n]$, converges to an element $\alpha$ of $K$. Moreover,  
\[
\displaystyle{\left| \alpha -\frac{s_{n}}{t_{n}} \right|=\frac{1}{|t_{n}||t_{n+1}|}}.
\]
\end{lemma}
\begin{proof}
Note that,
\[
\displaystyle{\frac{s_{n+1}}{t_{n+1}}-\frac{s_{n}}{t_{n}}=\frac{s_{n+1}t_{n}-s_{n}t_{n+1}}{t_{n}t_{n+1}}=\frac{(-1)^{n}}{t_{n+1}t_{n}}}\ \ \text{by (\ref{*})}.
\]
It is also easy to see that
\[
|t_n|=|c_nc_{n-1} \cdots c_1|, \ n \geq 1
\]
which in turn implies that $|t_n|$ is an increasing sequence as $c_j \in Z^{*}$ for $j \geq 1$. Then it follows that $(|t_{n+1}t_{n}|)$ is an increasing sequence as well. Now, using the properties of ultrametric absolute value, it can be easily seen that
\[
\left| \displaystyle{\frac{s_{m}}{t_{m}}}-\displaystyle{\frac{s_{n}}{t_{n}}} \right| = \frac{1}{|t_{n+1}t_{n}|}  \ \text{for any} \ m>n.
\]
As $(|t_{n+1}t_n|)$ is increasing, it follows that $(\frac{s_n}{t_n})$ is a Cauchy sequence, and hence converges to some $\alpha \in K$.
\end{proof}
Now, given any $\alpha \in K$, we generate its continued fraction expansion as follows:
\[
\alpha_0=\alpha, \ \alpha_{n+1}=(\alpha_{n}-\lfloor \alpha _n \rfloor )^{-1}, \ c_n=\lfloor \alpha _n \rfloor.
\]
If $\alpha_n =\lfloor \alpha _n \rfloor$ for some $n$, then $\alpha_{n+1}$ is not defined and the sequences $(\alpha_n)$ and $(c_n)$ are finite. Otherwise, two infinite sequences are generated by the above construction. Here, $c_n$'s the partial quotients and $\alpha_n$'s the complete quotients corresponding to the continued fraction expansion of $\alpha$. It is easy to see that $\alpha=[c_0;c_1,\ldots,c_n,\alpha_{n+1}]$.

Now, suppose $\alpha \in K$ be such that the sequences $(\alpha_n)$ and $(c_n)$ are infinite, and let
\[
\frac{s_n}{t_n}=[c_0;c_1,\ldots,c_n].
\]
We have $c_j \in Z^{*}$ for $j \geq 1$ by (\ref{**}). Then it follows from Lemma \ref{+} that the sequence of convergents $\frac{s_n}{t_n}$ converges to $\alpha$.
\begin{remark}
The definition of continued fraction in this article differs from the definition of continued fractions discussed by Capuano et al. in \cite{CMT2022}. Our defining conditions of floor function are less restrictive; in fact, we do not impose a condition like the $2$nd condition in Definition $3.1$ of \cite{CMT2022}. Also, our algorithm is less abstract which enables us to discuss the metrical theory of the associated continued fraction map. The following example shows that the $2$nd condition of Definition $3.1$ of \cite{CMT2022} may not be satisfied in our setup. 
\end{remark}

\begin{example}
Let $K=\Q_p(\beta)=\Q_5\left(\frac{1}{\sqrt{15}}\right)$. Then $[K:\Q_5]=2$ and $K$ is a totally ramified extension of $\Q_5$. If $\alpha\in K$ is given by $\alpha=\left(\sum \limits _{n=0}^{\infty} 5^{n}\right) \frac{1}{\sqrt{15}}$, then 
\[
 \lfloor \alpha \rfloor=\left \lfloor \left(\sum \limits _{n=0}^{\infty}5^{n}\right) \frac{1}{\sqrt{15}}\right \rfloor=\frac{1}{\sqrt{15}}.
\]
Now, let $|\ .\ |_{\mathfrak{v}*}$ be an ultrametric normalized absolute value  on the number field $\Q\left(\frac{1}{\sqrt{15}}\right)$ such that $\mathfrak{v}*$ is a non-Archimedean place lying over the prime $3$.
Then
\[
|\lfloor \alpha \rfloor|_{\mathfrak{v}*}=\left| \frac{1}{\sqrt{15}} \right|_{\mathfrak{v}*}=\left|\frac{1}{15}\right|^{\frac{1}{2}}_{3}=\sqrt{3}>1
\]
which violates the $2$nd condition of Definition $3.1$ of \cite{CMT2022}.
\end{example}

One of the main difficulties for continued fractions in the $p$-adic setup is that rational numbers do not necessarily have finite continued fraction expansions (also known as finiteness property) for many algorithms. In $1978$, Browkin modified Ruban's algorithm to achieve the finiteness property for $p$-adic continued fraction. The fact that Euclidean absolute value of the partial quotients in Browkin's algorithm is less than $\frac{p}{2}$ was crucially used in Browkin's proof of finiteness. In our setup, we prove the finiteness property for some small degree extensions of $\Q_p$ in the case of generalization of Browkin's algorithm, i.e., in the case $X=X_2$.

Let $p$ be either $3$ modulo $4$ or $5$ modulo $12$. In the first case, we take $K=\Q_p(\iota)$, where $\iota$ is the root of the polynomial $X^2+1=0$. In the $2$nd case we take $K=\Q_p(\omega)$, where $\omega$ is the root of the irreducible polynomial $X^2 + X + 1=0$. We show that the finiteness property holds in the cases of these extensions of $\Q_p$. Note that when $p \equiv -1 \pmod{12}$, then $\Q_p(\iota)$ and $\Q_p(\omega)$ gives rise to the same extension of $\Q_p$. For a cyclotomic extension $\Q_p(\gamma)$ of $\Q_p$, where $\gamma$ is some primitive $n$th root of unity, we define the Galois height of field rational elements as follows: for $\alpha \in \Q(\gamma)$,
\[
H(\alpha) := \max_{\sigma}|\sigma (\alpha)|_{\infty}
\]
where the maximum is taken over all (distinct modulo conjugation) Galois embeddings of $\Q(\gamma)$ inside $\C$, and $|y|_{\infty}$ denotes the Euclidean norm of the complex number $y$. The following lemma gives us the required bound on the Galois heights which will be useful in proving the finiteness property.
\begin{lemma} \label{lem 11}
If $b_0, b_1\ \in\ \Z \left[ \frac{1}{p} \right]\,\cap\,\left(-\frac{p}{2},\frac{p}{2}\right)$, then there exists $\delta >0$ such that $H\,(\,b_0+ b_1\iota\,)\ <\ p - \frac{1}{p} -\delta$ and $H\,(\,b_0 + b_1\omega\,)\ <\ p - \frac{1}{p}-\delta$.
\end{lemma}
\begin{proof}
Here, $\abs{b_j}_{\infty} < p/2$ for $j = 0, 1$ giving us that
\[
\abs{\sigma ( b_0 + b_1\iota )}_{\infty} = \abs{b_0 + b_1\sigma(\iota)}_{\infty} \leq \max \{ \abs{b_0}_{\infty}, \abs{b_1}_{\infty} \}\cdot\sqrt{2} < \frac{p}{\sqrt{2}} 
\]
for all $\sigma \in \operatorname{Gal}_{\Q} \big(\,\Q (\iota)\,\big)$ when $p \equiv 3 \pmod{4}$. Again,
\[
\abs{\sigma ( b_0 + b_1\omega )}_{\infty} = \abs{b_0 + b_1\sigma(\omega)}_{\infty} \leq \max \{ \abs{b_0}_{\infty}, \abs{b_1}_{\infty} \}\cdot\sqrt{3} <p \cdot \frac{\sqrt{3}}{2}
\]
for all $\sigma \in \operatorname{Gal}_{\Q} \big(\,\Q (\omega)\,\big)$ when $p \equiv 5 \pmod{12}$. Then it is clear that we can find a suitable $\delta>0$ such that the assertions of the lemma hold.

\end{proof}
\begin{proposition}\label{prop 12}
Let $K=\Q_p(\gamma)$ where $\gamma=\iota$ when $p \equiv 3 \pmod{4}$ and $\gamma=\omega$ when $p \equiv 5 \pmod{12}$. Also let $X=X_2$, i.e., the partial quotients of the continued fraction expansion of any element of $\Q_p(\gamma)$ are elements of the form
\[
b_0+b_1\gamma \ \text{with} \ b_0,b_1 \in \Z \left[\frac{1}{p} \right] \cap \left( -\frac{p}{2}, \frac{p}{2} \right).
\]
Then, any $\alpha \in \Q(\gamma)$ has a finite continued fraction expansion.
\end{proposition}
\begin{proof}
We use a suitable generalisation of the method used in Proposition $4.3$ of \cite{CMT2022}. As $\alpha \in \Q(\gamma)$, we can express $\alpha$ as
\[
\alpha=\frac{X_{0}}{Y_{0}}
\]
with $X_0 \in \Z \left[\frac{1}{p}\right][\gamma]$, $Y_0 \in \Z$ and $p \not|\  Y_0$. We define two sequences $(U_n)$ and $(Y_n)$ as follows:
\[
U_n=s_n-\alpha t_n, \ Y_n=Y_0U_n.
\]
Then, it is easy to see that 
\[
U_n=(-1)^{n+1}\prod \limits_{j=1}^{n+1}\frac{1}{\alpha_j}
\]
and, consequently,
\begin{equation}\label{eq*}
|U_n|=\prod \limits_{j=1}^{n+1}\frac{1}{|c_j|} \ \text{as} \ |\alpha_j|=|c_j|.
\end{equation}
It is also easy to see that the sequence $(Y_n)$ satisfies the recurrence relation
\begin{equation}\label{eq**}
Y_n=c_n Y_{n-1}+Y_{n-2}.
\end{equation}
Clearly, $Y_n \in \Z \left[\frac{1}{p}\right][\gamma]$. Also, by definition of $Y_n$ and (\ref{eq*}), we have $Y_n \in p^n \Z _p[\gamma]$. Hence, $Y_n \in \Z \left[\frac{1}{p}\right][\gamma] \cap p^n \Z _p[\gamma]=p^n \Z [\gamma]$. Taking the Galois height of both sides of (\ref{eq**}), and then dividing by $p^n$, we have
\begin{align*}
\frac{H(Y_n)}{p^n}& \leq \frac{H(c_n)}{p}\frac{H(Y_{n-1})}{p^{n-1}}+\frac{1}{p^2}\frac{H(Y_{n-2})}{p^{n-2}}\\&<\left( p-\frac{1}{p}-\delta \right) \cdot \frac{1}{p} \cdot \frac{H(Y_{n-1})}{p^{n-1}}+\frac{1}{p^2} \cdot \frac{H(Y_{n-2})}{p^{n-2}}.
\end{align*}
Let $T_n=\frac{H(Y_n)}{p^n}$, $D_1=\left( p-\frac{1}{p}-\delta \right) \cdot \frac{1}{p}$, $D_2=\frac{1}{p^2}$. Then 
\[
T_n<D_1\ T_{n-1}+D_2\ T_{n-2}.
\]
Since $D_0+D_1<1$, it follows from Lemma $4.2$ of \cite{CMT2022} that $|T_n|_{\infty} \to 0$ as $n \to \infty$. Hence there exists $n_0\in\N$ such that $Y_n=0$ $\forall \, n \geq n_0$ since $\frac{Y_n}{p^n} \in \Z[\gamma]$. This means that $\alpha=\frac{s_n}{t_n}$ for some $n$, and consequently, $\alpha$ has a finite continued fraction expansion.
\end{proof}

  \section{Exactness}
Let $K=\Q_p(\gamma, \beta)$ be a finite extension of $\Q_p$, and $\lfloor\ .\ \rfloor$ be the floor function defined in \ref{E:floor}. The continued fraction map $T$ is defined on $B(0,1)$ inside $K$, as follows:
    \begin{equation}\label{E:cfm}
      T(\alpha)=\frac{1}{\alpha}- \left \lfloor \frac{1}{\alpha} \right \rfloor \text{ for } \alpha \neq 0 \text{ and } T(0) = 0,
    \end{equation}
    where $\lfloor \cdot \rfloor$ is as def{i}ned in~\eqref{E:floor}. 
In this section, we shall prove the exactness of $T$, and in the subsequent section we prove various metrical results as consequences of exactness. We shall be considering the continued fraction map corresponding to the extension of Ruban's algorithm in finite extensions of $\Q_p$, though similar assumptions hold for the continued fraction map corresponding to the extension of the Browkin's algorithm as well. Now, let $\alpha\in B(0,1)$ and $\alpha = [0; c_1,c_2,\ldots]$ be the continued fraction axpansion of $\alpha$. To emphasize the dependence on $\alpha$, we will also use $c_k(\alpha)$ to denote the $k$th partial quotient of the continued fraction expansion of $\alpha$, i.e., $\alpha = [0;c_1(\alpha),c_2(\alpha),\ldots]$. Note that,
    \[
      T^{n}(\alpha)=[ 0; c_{n+1}(\alpha),c_{n+2}(\alpha),\ldots], \ \text{and} \  c_k \big( T^{n}(\alpha) \big) = c_{n+k}(\alpha)
    \]
    for all $k \geq 1$ and $n \geq 0$.

 Recall that a measure preserving dynamical system $(X, \mathcal{C}, \nu, S)$ is said to be exact if
 \[
  \bigcap \limits_{n=0}^{\infty}S^{-n} \mathcal{C}=\mathcal{N} \ (\text{mod} \ \nu),
  \]
 where $\mathcal{N}$ is the trivial sub $\sigma$-algebra of $\mathcal{C}$ generated by the sets of measure $0$ or $1$.
 For $n \in \N$, and $c_1,\ldots,c_n \in Z^*$, let $\Delta_{c_1,\ldots,c_n}$ denote the \emph{cylinder set of length $n$}, i.e., 
 \begin{equation}\label{E:cyl} \Delta_{c_1,\ldots,c_n} = \big\{ [0; c_1,\ldots,c_{n-1},c_n+\beta ]\ :\ \beta \in B(0,1) \big\}.
\end{equation}
    The following lemma gives an alternate description of a cylinder set which will be helpful in calculating its measure. The proof of this lemma is similar to the proof of Lemma $2$ of \cite{LN14}.
    \begin{lemma}\label{L:LN14}
      For any f{i}nite sequence $c_1,\ldots,c_n \in Z^*$,
      \[
        \Delta_{c_1,\ldots,c_n} = B\,\big(\,[0; c_1,\ldots,c_n],\,|c_1 \cdots c_n|^{-2}\,\big).
      \]
    \end{lemma}

\noindent Because of the above lemma, it is not hard to see that the Borel $\sigma$-algebra on $B(0,1)$ is generated by the cylinder sets described above. We denote by $\mathcal{B}$ the Borel $\sigma$-algebra on $B(0,1)$. Also, let $\mu$ be the restriction of the Haar measure on $B(0,1)$, and $T$ be the continued fraction map on $B(0,1)$ defined above. We first show that $T$ is measure-preserving. Note that two cylinders $\Delta_{c_1,\ldots,c_n}$ and $\Delta_{d_1,\ldots,d_n}$ of the same length are disjoint if and only if $c_j \neq d_j$ for some $1 \leq j \leq n$.
  \begin{lemma}\label{E:mpds}
    The dynamical system $\big(B(0,1),\mathcal{B},\mu,T\big)$ is measure-preserving.
  \end{lemma}
  \begin{proof}
    Since the cylinder sets generate the Borel $\sigma$-algebra, it is enough to show that $T$ is measure-preserving on cylinder sets. For any cylinder set $\Delta_{c_1,\ldots,c_n}$, there exists $s \in \mathbb{Z}$ and $i \in \{ 0,\ldots,e-1 \} $ such that $|c_1 \cdots c_n|^{-2}=p^{s\,+\,\frac{i}{e}}$, and consequently, 
    \[
    \mu(\Delta_{c_1,\ldots,c_n})=\mu (B([0,c_1,\ldots,c_n],|c_1 \cdots c_n|^{-2}))=p^{ms+fi}.
    \]
    The inverse image of $\Delta_{c_1,\ldots,c_n}$ under $T$ is given by a disjoint union as follows:
    \begin{equation}\label{E:decom}
        T^{-1}\Delta_{c_1, \ldots, c_n}=\bigcup_{c \in Z^*}\Delta_{c,c_1, \ldots, c_n}.
    \end{equation}
    Then
    \begin{align*}
    \mu(T^{-1}\Delta_{c_1,\ldots,c_n})&=\sum \limits_{c \in Z^*}\mu (B([c,c_1,\ldots,c_n],|c|^{-2}p^{s+\frac{i}{e}}))\\&=\sum \limits_{n=1}^{\infty}p^{fn}(p^f-1)p^{ms+fi-2fn} \ (\text{using} \ (\ref{Z*:coun}))\\&=(p^f-1)\frac{p^{ms+fi-f}}{1-\frac{1}{p^f}}\\&=p^{ms+fi}\\&=\big(p^{s+\frac{i}{e}}\big)^{m}\\&=\mu(\Delta_{c_1,\ldots,c_n}).
    \end{align*}
\end{proof}
The following technical lemma which is analogous to Lemma $4$ of \cite{LN14}, is a crucial ingredient in proving exactness of the continued fraction map.  
\begin{lemma}\label{E:imm}
    For the dynamical system $\big(B(0,1),\mathcal{B},\mu,T\big)$, if $E \in \mathcal{B}$, then for any natural number $n$ and cylinder set $\Delta_{c_1,\ldots,c_n}$, we have 
    \[
    \mu(\Delta_{c_1,\ldots,c_n} \cap T^{-n}E)=\mu (\Delta_{c_1,\ldots,c_n})\mu(E).
    \]
\end{lemma}
\begin{proof}
    It is enough to consider $E$ to be a cylinder set. Let $E=\Delta_{d_1,\ldots,d_m}$. Then there exist $s_{1} ,s_{2} \in \Z$ and $i_{1},i_{2} \in \{ 0,\ldots,e-1 \}$ such that $|c_1 \cdots c_n|^{-2}=p^{s_{1}+\frac{i_{1}}{e}}$ and $|d_1 \cdots d_m|^{-2}=p^{s_{2}+\frac{i_{2}}{e}}$. Now,
    \[
    T^{-n}\Delta_{d_1,\ldots,d_m}=\bigcup_{c'_1,\ldots,c'_n \in Z^*}\Delta_{c'_1,\ldots,c'_n,d_1,\ldots,d_m}.
    \]
    Also, $\Delta_{c_1,\ldots,c_n} \cap T^{-n}\Delta_{d_1,\ldots,d_m}=\mu(\Delta_{c_1,\ldots,c_n,d_1,\ldots,d_m})$. Then,
    \begin{align*}
        \mu ( \Delta_{c_1,\ldots,c_n} \cap T^{-n} \Delta_{d_1,\ldots,d_m})&=\mu (\Delta_{c_1,\ldots,c_n,d_1,\ldots,d_m})\\&=\mu(B([c_1,\ldots,c_n,d_1,\ldots,d_m],|c_1 \cdots c_n d_1 \cdots d_m|^{-2}))\\&=\mu(B([c_1,\ldots,c_n,d_1,\ldots,d_m],p^{s_{1}+\frac{i_{1}}{e}}p^{s_{2}+\frac{i_{2}}{e}}))\\&=p^{ms_1+fi_1} \cdot p^{ms_2+fi_2}\\&=\mu(\Delta_{c_1,\ldots,c_n}) \mu( \Delta_{d_1,\ldots,d_m}).
    \end{align*}
\end{proof}
Now we show that the continued fraction map $T$ is exact.
\begin{theorem}\label{E:eds}
The dynamical system $\big(B(0,1),\mathcal{B},\mu,T\big)$ is an exact dynamical system.
\end{theorem}
\begin{proof}
It is enough to show that $\bigcap_{n=0}^{\infty}T^{-n}\mathcal{B} \subseteq \mathcal{N}$. Let $E \in \bigcap_{n=0}^{\infty}T^{-n}\mathcal{B}$. Then for each $n \geq 1$, there exists $E_n \in \mathcal{B}$ such that $E = T^{-n}E_n$ and $\mu(E_n)=\mu(E)$. Now, for each cylinder set $\Delta_{c_1,\ldots,c_n}$ of length $n$, 
\begin{align*}
\mu(E \cap \Delta_{c_1,\ldots,c_n})&=\mu(T^{-n}E_n \cap \Delta_{c_1,\ldots,c_n})\\&=\mu(E)\mu(\Delta_{c_1,\ldots,c_n}) \ (\text{by Lemma } \ \ref{E:imm}).
\end{align*}
Then it follows from Lemma $5$ of \cite{LN14} that $\mu(E)=0$ or $1$, consequently, $E\in \mathcal{N}$. 
\end{proof}

\section{Metrical results}
Now we obtain results analogous to the metrical results of \cite{LN14} in our setup. Since $T$ is exact, it is weak-mixing as well, i.e.,
\[
\frac{1}{n}\sum \limits_{k=1}^{n}\big| \mu (E \cap T^{-k}F)-\mu (E) \mu (F) \big| \to 0
\]
as $n \to \infty$ for any $E,F \in \mathcal{B}$.
Weak-mixing property of the continued fraction map enables one to consider metrical results in the context of certain subsequences. This is done in \cite{NR2} for continued fraction map in the case of real numbers, and in \cite{LN14} in the positive characteristic setup. We do a similar study here for continued fraction map on $B(0,1)$ inside $K$. Before proceeding further we recall two definitions which plays crucial role in the discussion of metrical theory using subsequences.
\begin{definition}
A strictly increasing sequence of positive integers $(a_n)_{n=1}^{\infty}$ is said to be $L^{2}$-\emph{good universal}, if for each dynamical system $(X, \mathcal{C}, \nu ,S)$ and $g \in L^{2}(X, \mathcal{C}, \nu)$, the limit 
\[
\lim \limits_{n \to \infty} \frac{1}{n} \sum \limits_{j=1}^{n}g(S^{a_j-1} \alpha)
\]
exists $\nu$-almost everywhere.    
\end{definition} 
\begin{definition}
A sequence of real numbers $(x_n)_{n=1}^{\infty}$ is called \emph{uniformly distributed modulo 1}, if for each interval $I \subseteq [0,1)$, we have 
 \[
\lim \limits_{n \to \infty}\frac{1}{n} \cdot \# \{ 1 \leq j \leq n : \{ x_j \} \in I \} =|I|,
\]
where $|I|$ denotes the length of $I$ and $\{ x_j \}$ denotes the fractional part of $x_j$.    
\end{definition}
Please see \cite{LN14} for examples of $L^2$-good universal sequences. The following proposition is a consequence of weak-mixing, the proof of which can be found in \cite{NR2}.

\begin{proposition}\label{Prop*}
Let $(X, \mathcal{C}, \nu, S)$ be a weak-mixing dynamical system. Suppose $(a_n)^{\infty}_{n=1}$ is an $L^2$-good universal sequence of natural numbers such that $(a_n \gamma)^{\infty}_{n=1}$ is uniformly distributed modulo $1$ for any irrational number $\gamma$. Then for any $g \in L^2(X, \mathcal{C}, \nu)$,
\[
\lim \limits_{n \to \infty} \frac{1}{n} \sum \limits_{j=1}^{n}g\big(S^{a_{j}-1}\alpha\big)=\int \limits_{X} g\ \d \nu
\]
$\nu$-almost everywhere.
\end{proposition}
\begin{proposition}\label{Prop**}
Let $F: \mathbb{R}_{\geq 0} \to \mathbb{R}$ be an increasing function such that 
\[
\int \limits_{B(0,1)} \big|F(|c_1(\alpha)|)\big|^2 \, \d \mu < \infty .
\]
For any natural number $n$ and non-negative real numbers $d_1,\ldots,d_n$, let the generalized average be defined as 
\[
M_{F,n}(d_1,\ldots,d_n)=F^{-1} \left( \frac{F(d_1)+\cdots +F(d_n)}{n} \right).
\]
If $(a_n)^{\infty}_{n=1}$ is an $L^2$-good universal sequence of natural numbers such that $(a_n \gamma)^{\infty}_{n=1}$ is uniformly distributed modulo $1$ for any irrational number $\gamma$, then 
\[
\lim \limits_{n \to \infty} M_{F,n}(|c_{a_{1}}(\alpha)|,\ldots,|c_{a_{n}}(\alpha)|)=F^{-1} \left( \int \limits_{B(0,1)} F(|c_1(\alpha)|) \, \d \mu \right)
\]
$\mu$-almost everywhere.
\end{proposition}
\begin{proof}
Apply Proposition \ref{Prop*} to the function $g(\alpha)=F\big(|c_{1}(\alpha)|\big)$.
\end{proof}
The following is also a consequence of Proposition \ref{Prop*}, in which one considers a function from $\mathbb{R}^{k}_{\geq 0}$ to $\mathbb{R}$.
\begin{proposition}\label{Prop***}
Suppose that $H:\mathbb{R}^{k}_{\geq 0} \to \mathbb{R}$ is a function such that 
\[
\int \limits_{B(0,1)} \Big|H\big(|c_1(\alpha)|,\ldots,|c_k(\alpha)|\big)\Big|^2 \, \d \mu < \infty ,
\]
and if $(a_n)^{\infty}_{n=1}$ is an $L^2$-good universal sequence of natural numbers such that $(a_n \gamma)^{\infty}_{n=1}$ is uniformly distributed modulo $1$ for any irrational number $\gamma$. Then 
\[
\lim \limits_{n \to \infty} \frac{1}{n} \sum \limits_{j=1}^{n}H\big(|c_{a_{j}}(\alpha)|,\ldots,|c_{a_{j}+k-1}(\alpha)|\big)=\sum \limits_{(i_1,\ldots,i_k) \in \mathbb{N}^k} H\left(p^{\frac{i_1}{e}},\ldots,p^{\frac{i_k}{e}}\right) \frac{(p^f-1)^k}{p^{f(i_1+\cdots +i_k)}}
\]
$\mu$-almost everywhere.
\end{proposition}
\begin{proof}
Let $g(\alpha)=H\big(|c_1(\alpha)|,\ldots,|c_k(\alpha)|\big)$. Then applying Proposition \ref{Prop*} to the function $g$, we get
\begin{align*}
\lim \limits_{n to \infty} \frac{1}{n} \sum \limits_{j=1}^{n}H\big(|c_{a_{j}}(\alpha)|,\ldots,|c_{a_{j}+k-1}|\big)&=\int \limits_{B(0,1)} H\big(|c_1(\alpha)|,\ldots,|c_{k}(\alpha)|\big) \,\d \mu \\&=\sum \limits_{(i_1,\ldots,i_k) \in \mathbb{N}^k} H\left(p^{\frac{i_1}{e}},\ldots,p^{\frac{i_k}{e}}\right)\frac{(p^f-1)^k}{p^{f(i_1+\cdots +i_k)}},
\end{align*}
where we have used (\ref{Z*:coun}) to get the last equality.
\end{proof}
Now, we calculate the asymptotic frequency of partial quotients being some particular element of $Z^{*}$.
\begin{lemma}
If $(a_n)$ is a sequence as in the above propositions, then for any $z \in Z^{*}$, 
\[
\lim \limits_{n \to \infty} \frac{1}{n} \cdot \# \big\{ 1 \leq j \leq n :c_{a_{j}}(\alpha)=z \big\} =\frac{1}{|z|^{2m}}
\]
almost everywhere with respect to $\mu$.
\end{lemma}
\begin{proof}
Applying Proposition \ref{Prop*} with $g(\alpha)= \chi _{ \{ z \} }\big(c_1(\alpha)\big)$, we have
\begin{align*}
\lim \limits_{n \to \infty} \frac{1}{n} \cdot \# \big\{ 1 \leq j \leq n :c_{a_{j}}(\alpha)=z \big\} &=\int \limits_{B(0,1)} \chi _{ \{ z \} }\big(c_1(\alpha)\big) \, \d \mu \\&= \mu \big( \{ \alpha \in B(0,1) : c_{1}(\alpha)=z \} \big)\\&= \mu \left( B \left( \frac{1}{z},|z|^{-2} \right) \right) \\&= \frac{1}{|z|^{2m}}.
\end{align*}
\end{proof}
In the following two results, we assume that $(a_n)^{\infty}_{n=1}$ is an $L^2$-good universal sequence of natural numbers such that $(a_n \gamma)^{\infty}_{n=1}$ is uniformly distributed modulo $1$ for any irrational number $\gamma$. The next result is a version of Khinchin's theorem regarding the geometric mean of the partial quotients in the case of real continued fraction.
\begin{proposition}
For almost every $\alpha \in B(0,1)$ with respect to the Haar measure,
\[
\lim \limits_{n \to \infty} \frac{1}{n} \sum \limits_{j=1}^{n} \big( -\mathfrak{v} \big(c_{a_j}(\alpha)\big)\big)=\frac{p^{f}}{(p^{f}-1)e}.
\]
\end{proposition}
\begin{proof}
Applying Proposition \ref{Prop*} to the function $g(\alpha)=\log _{p}\big(|c_1(\alpha)|\big)$, we get
\[
\lim \limits_{n \to \infty} \frac{1}{n} \sum \limits_{j=1}^{n} \log _p \big(|c_{a_j}|\big)
    =\int \limits_{B(0,1)} \log _p \big(|c_1 (\alpha)|\big) \, \d \mu . 
\]
Now,
\[
\sum \limits_{j=1}^{n} \log _p \big(|c_{a_j}|\big)=\sum \limits_{j=1}^{n} \log _p \Big(p^{-\mathfrak{v}\big(c_{a_{j}}(\alpha)\big)}\Big)=\sum \limits_{j=1}^{n}-\mathfrak{v}\big(c_{a_{j}}(\alpha)\big).
\]
Also,
\begin{align*}
\int \limits_{B(0,1)} \log _p\big(|c_1(\alpha)|\big) \, \d \mu &=\sum \limits_{n=1}^{\infty} \frac{n}{e} \mu \left\{ \alpha \in B(0,1) :|\alpha|=p^{-\frac{n}{e}} \right\} \\&= \sum \limits_{n=1}^{\infty} \frac{n}{e} \Big( p^{\big(-\frac{n}{e}+\frac{1}{e}\big)m}-p^{-\frac{n}{e}\cdot m} \Big) \\&= \sum \limits_{n=1}^{\infty}p^{-nf}\big(p^f-1\big)\\&=\frac{p^f-1}{e} \cdot \sum \limits_{n=1}^{\infty}n \cdot p^{-nf} \\&=\frac{p^f-1}{e} \cdot \frac{p^{-f}}{\big(1-p^{-f}\big)^{-2}}\\&=\frac{p^f}{\big(p^f-1\big)e}
\end{align*}
Then the proposition follows.
\end{proof}

In the following theorem, we find the asymptotic frequency of partial quotients taking some specified absolute value (or greater or equal to some specified absolute value or absolute values in certain range).
\begin{theorem}
For any positive integer $l$,\\
$(i)$ $\lim \limits_{n \to \infty} \frac{1}{n} \cdot \# \left \{ 1 \leq j \leq n : |c_{a_{j}}| = p^{\frac{l}{e}} \right \} =\frac{\displaystyle p^{f}-1}{\displaystyle{p^{f \, l}}}$\\
and $(ii)$ $\lim \limits_{n \to \infty} \frac{1}{n} \cdot \# \left \{ 1 \leq j \leq n : |c_{a_{j}}| \geq p^{\frac{l}{e}} \right \} =\frac{\displaystyle 1}{\displaystyle{p^{f(l-1)}}}$ $\mu$-almost everywhere.\\
If $k$ is another positive integer with $k<l$, then \\
$(iii)$ $\lim \limits_{n \to \infty} \frac{1}{n} \cdot \# \left \{ 1 \leq j \leq n : p^{\frac{k}{e}} \leq |c_{a_{j}}| < p^{\frac{l}{e}} \right \} =\frac{\displaystyle 1}{\displaystyle{p^{f(k-1)}}} \left( 1- \frac{\displaystyle 1}{\displaystyle{p^{f(l-k)}}} \right)$\\
$\mu$-almost everywhere.\\
\end{theorem}
\begin{proof}
Apply Proposition \ref{Prop*} with $g(\alpha)=\chi _{ \big\{ p^{\frac{l}{e}} \big\} }\big(|c_1(\alpha)|\big)$ for the proof of $(i)$. Similarly, for the proof of $(ii)$ and $(iii)$, apply Proposition \ref{Prop*} with $g(\alpha)=\chi _{ \big[p^{\frac{l}{e}}, \infty \big)}\big(|c_1(\alpha)|\big)$ and $g(\alpha)=\chi _{ \big[p^{\frac{k}{e}},p^{\frac{l}{e}}\big)}\big(|c_1(\alpha)|\big)$, respectively.
\end{proof}

Given a measure-preserving transformation $S$ on a probability measure space $(X, \mathcal{C}, \nu)$, we know that $\frac{1}{n}\sum \limits_{j=0}^{n-1}g(S^jx)$ converges almost everywhere. But what happens if we consider moving averages? This means, given a sequence of pairs of positive integers $(a_n,b_n)^{\infty}_{n=1}$ what can be said about the convergence of $\frac{1}{b_n}\sum \limits_{j=0}^{b_n-1}g(S^{a_{n}+j}x)$ for almost every $x$. In \cite{ARJ90}, necessary and sufficient conditions were given for this kind of moving averages to converge.

Let $\Omega$ be an infinite collection of points inside $\mathbb{Z} \times \mathbb{N}$. Define
\begin{align*}
&\Omega^h= \{ (n,k):(n,k) \in \Omega \ \text{and} \ k \geq h \} \\&  \Omega^{h}_{\alpha}= \{ (z,s) \in \Z ^2 : |z-y|<\alpha \, (s-r) \ \text{for some} \ (y,r) \in \Omega^h \} \\& \Omega^{h}_{\alpha}(\lambda)= \{ n: (n,\lambda) \in \Omega^{h}_{\alpha} \} .
\end{align*}
Following \cite{LN14}, we call a sequence of pairs of natural numbers $(a_n,b_n)$ Stoltz if there exists a function $h=h(t)$ tending to infinity with $t$ such that
\[
(a_n,b_n)^{\infty}_{n=t} \in \Omega^{h(t)}
\]
for some $\Omega$ inside $\mathbb{Z} \times \mathbb{N}$, and $\exists$ $l, \alpha$ and $A_{\alpha}>0$ such that 
\[
\left|\Omega^{l}_{\alpha}(\lambda)\right| \leq A_{\alpha}(\lambda),
\]
where $\left|\Omega^{l}_{\alpha}(\lambda)\right|$ denotes the cardinality of the set $\Omega^{l}_{\alpha}(\lambda)$. 

The following proposition which can be considered as the moving average version of Proposition \ref{Prop*}, is the base of the metrical results corresponding to the continued fraction map $T$ in the context of moving averages. The proof of the proposition is essentially contained in \cite{ARJ90}.
\begin{proposition}
Let $(X, \mathcal{C}, \nu, S)$ be an ergodic dynamical system, and let $(a_n,b_n)^{\infty}_{n=1}$ be a Stoltz sequence of natural numbers. Then for any $g \in L^1(X, \mathcal{C}, \nu)$,
\[
\lim \limits_{n \to \infty} \frac{1}{b_n} \sum \limits_{j=1}^{b_n}g\big(S^{a_n+j-1}\alpha\big)=\int \limits_{X} g \, \d \nu
\]
$\nu$-almost everywhere.
\end{proposition}
\noindent The readers are referred to \cite{LN14} for some examples of Stoltz sequence, and to \cite{ARJ90} for some sequence of pairs of natural numbers for which the assumption of the above proposition fails. 

Now we state the metrical results using moving averages. The results are analogous to the results mentioned above, and proofs are similar.
\begin{proposition}
Let $F: \R _{\geq 0} \to \R$ and $M_{F,n}$ be as in Proposition \ref{Prop**}, and let $(a_n,b_n)$ be a Stoltz sequence of pairs of natural numbers. Then
\[
\lim \limits_{n \to \infty}M_{F,n}\big(|c_{a_{n}+1}(\alpha)|,\ldots,|c_{a_{n}+b_{n}}(\alpha)|\big)=F^{-1} \left(\ \int \limits_{B(0,1)} F\big(|c_1 (\alpha)|\big) \ \d \mu\ \right)
\]
almost everywhere with respect to $\mu$.
\end{proposition}
\begin{proposition}
Suppose $H: \R^{k} _{\geq 0} \to \R$ is a function as in Proposition \ref{Prop***}, and $(a_n,b_n)$ be a Stoltz sequence of pairs of natural numbers. Then
\[
\lim \limits_{n \to \infty} \frac{1}{b_n} \sum \limits_{j=1}^{b_n}H\big(|c_{a_{n}+j}(\alpha)|,\ldots,|c_{a_{n}+j+k-1}(\alpha)|\big)=\sum \limits_{(i_1,\ldots,i_k) \in \mathbb{N}^k} H\left(p^{\frac{i_1}{e}},\ldots,p^{\frac{i_k}{e}}\right) \frac{\big(p^f-1\big)^k}{p^{f(i_1+\cdots +i_k)}}
\]
$\mu$-almost everywhere.
\end{proposition}
\noindent We include the other results in the following theorem.
\begin{theorem}
Let $(a_n,b_n)$ be a Stoltz sequence of pairs of natural numbers. We consider all the statements mentioned below in the almost everywhere sense with respect to the measure $\mu$.
\begin{align*}
&(i) \ \text{For any} \ z \in Z^{*},\\& \ \ \ \lim \limits_{n \to \infty} \frac{1}{b_n} \ \# \big\{ 1 \leq j \leq b_n :c_{a_{n}+j}(\alpha)=z \big\} =\frac{1}{|z|^{2m}}, \\& (ii) \lim \limits_{n \to \infty} \frac{1}{b_n} \sum \limits_{j=1}^{b_n}  -\mathfrak{v} \big(c_{a_{n}+j}(\alpha)\big)=\frac{p^{f}}{\big(p^{f}-1\big)e}. \\& (iii) \ \text{For any} \ l \in \N, \\& \ \ \ \ \ \lim \limits_{n \to \infty} \frac{1}{b_n} \ \# \left \{ 1 \leq j \leq b_n : |c_{a_{n}+j}| = p^{\frac{l}{e}} \right \} =\frac{\displaystyle p^{f}-1}{\displaystyle{p^{f \, l}}}, \\& \ \ \ \ \ \ \text{and} \ \lim \limits_{n \to \infty} \frac{1}{b_n} \ \# \left \{ 1 \leq j \leq b_n : |c_{a_{n}+j}| \geq p^{\frac{l}{e}} \right \} =\frac{\displaystyle 1}{\displaystyle{p^{f(l-1)}}}. \\& (iiii) \ \text{For} \ k,l \in \N \ \text{with} \ k<l,\\& \ \ \ \ \ \ \lim \limits_{n \to \infty} \frac{1}{b_n} \ \# \left \{ 1 \leq j \leq b_n : p^{\frac{k}{e}} \leq |c_{a_{n}+j}| < p^{\frac{l}{e}} \right \} =\frac{\displaystyle 1}{\displaystyle{p^{f(k-1)}}} \left( 1- \frac{\displaystyle 1}{\displaystyle{p^{f(l-k)}}} \right).
\end{align*}
\end{theorem}
\begin{acknowledgement}
Both the authors are thankful to Dr. L. Singhal for many helpful discussions and some insightful suggestions especially in  Lemma \ref{L:1ud}. He also suggested using the bounds in Lemma \ref{lem 11} in the proof of Proposition \ref{prop 12}.  Prashant J. Makadiya acknowledges the support of the Government of Gujarat through the SHODH (Scheme of Developing High-Quality Research) fellowship. Prashant J. Makadiya also thanks the Council of Scientific and Industrial Research (CSIR), India for their support through the CSIR-JRF fellowship.
\end{acknowledgement}

\bibliography{ref}
\bibliographystyle{plain}
\end{document}